\documentclass[12pt,draft]{amsart}

\usepackage{amscd}
\usepackage{pb-diagram}

\newtheorem{thm}{Theorem}[section]   
\newtheorem{cor}[thm]{Corollary}     
\newtheorem{lem}[thm]{Lemma}         
\newtheorem{prop}[thm]{Proposition}  

\theoremstyle{definition}

\newtheorem{defn}[thm]{Definition}   

\theoremstyle{remark}

\numberwithin{equation}{section}     

\newcommand{\thmref}[1]{Theorem~\ref{#1}}
\newcommand{\corref}[1]{Corollary~\ref{#1}}
\newcommand{\lemref}[1]{Lemma~\ref{#1}}
\newcommand{\propref}[1]{Proposition~\ref{#1}}

\newcommand{\ind}{\operatorname{Ind}}
\newcommand{\Ind}{\ind}

\newcommand{\rep}{\operatorname{Rep}}
\newcommand{\Rep}{\rep}

\newcommand{\supp}{\operatorname{\rm supp}}

\newcommand{\im}{\operatorname{\rm im}}
\newcommand{\ad}{\operatorname{\rm Ad}}

\newcommand{\aut}{\operatorname{Aut}}
\newcommand{\dashind}{\operatorname{\!-Ind}}
\renewcommand{\inf}{\operatorname{inf}}

\newcommand{\half}{\frac{1}{2}}

\renewcommand{\H}{\mathcal H}

\newcommand{\K}{\mathcal K}
\newcommand{\D}{\mathcal D}
\renewcommand{\L}{\mathcal L}

\newcommand{\cstar}{\ensuremath{C^*}-}
\renewcommand{\star}{\text{${}^*$}-}
\newcommand{\spacetext}[1]{\quad\text{#1}\quad}
\newcommand{\righttext}[1]{\qquad\text{#1 }}

\newcommand{\deltahat}{{\hat{\delta}}}
\renewcommand{\phi}{\varphi}

\newcommand{\lip}[3]{
  {\vphantom\langle}_{#1}\!\!\left\langle{#2},{#3}\right\rangle  }
\newcommand{\rip}[3]{
  \left\langle{#2},{#3}\right\rangle_{\!{#1}}  }

\begin{document}

\title{Imprimitivity for $C^*$-Coactions of Non-Amenable Groups}

\author{S. Kaliszewski}
\address{Department of Mathematics\\University of Newcastle\\
 Newcastle, New South Wales 2308\\Australia}
\email{kaz@frey.newcastle.edu.au}

\author{John Quigg}
\address{Department of Mathematics\\Arizona State University\\
 Tempe, Arizona 85287}
\curraddr{Department of Mathematics\\University of Newcastle\\
 Newcastle, New South Wales 2308\\Australia}
\email{quigg@math.la.asu.edu}

\thanks{This research was partially supported by the National Science
Foundation under Grant No. DMS9401253, and by the Australian Research
Council.}



\date{March 22, 1996 (revised)}

\begin{abstract}
We give a condition on a full coaction $(A,G,\delta)$
of a (possibly) nonamenable group $G$ and a closed normal
subgroup $N$ of $G$ which ensures that Mansfield imprimitivity works;
i.e. that $A\times_{\delta{\vert}} G/N$ is Morita equivalent to
$A\times_\delta G\times_{\deltahat,r} N$.  This condition obtains if
$N$ is amenable or $\delta$ is normal.  It is preserved under Morita
equivalence, inflation of coactions, the stabilization trick of
Echterhoff and Raeburn, and on passing to twisted coactions.
\end{abstract}

\maketitle


\section{Introduction}

For a nondegenerate
reduced coaction $\delta\colon A\to M(A\otimes C^*_r(G))$ and a
closed normal amenable subgroup $N$ of $G$, Mansfield
\cite[Theorem~27]{ManJF}\ showed that the iterated crossed product
$A\times_{\delta}G\times_{\hat\delta}N$ is Morita equivalent to the
restricted crossed product $A\times_{\delta{\vert}}G/N$.  
Later, Phillips and Raeburn \cite[Theorem~4.1]{PR-TC}\
proved that for a closed normal amenable subgroup
$N$ of $G$ and a twisted coaction of $(G,G/N)$ on $A$, there is an
action $\tilde\delta$ of $N$ on
the twisted crossed product $A\times_{G/N}G$ such that
$A\times_{G/N}G\times_{\tilde\delta}N$ is Morita equivalent to $A$; this
can be interpreted as a twisted version of Mansfield's theorem.

Both the above results are coaction versions of Green's imprimitivity
theorem for actions, and hence play an important role in the study of
induced representations of crossed products by coactions.  Naturally,
we would like to remove the amenability hypothesis on $N$, which is
only there to ensure that the restricted coaction $\delta{\vert}\colon A\to
M(A\otimes C^*_r(G/N))$ is well-defined.

Accordingly, in this paper we consider a nondegenerate full coaction
$\delta\colon A\to M(A\otimes C^*(G))$ and a closed normal subgroup $N$
of $G$.  We also use the abstract
characterization of the crossed product $A\times_\delta G$ using
universal properties, whereas Mansfield represents the crossed
product on Hilbert space.
In this situation, the proof of Mansfield's theorem shows that
$A\times G\times_r N$ is Morita equivalent to $\im(j_A\times j_G{\vert})$,
where $(j_A,j_G)$ is the canonical covariant homomorphism of
$(A,G,\delta)$ into $M(A\times G)$.  We conclude
(\corref{m-imprim-cor}) that if
$j_A\times j_G{\vert}\colon A\times G/N\to M(A\times G)$ is faithful, then
there is an $A\times G\times_r N$ -- $A\times G/N$ imprimitivity
bimodule $Y_{G/N}^G$, which we describe explicitly.  When $\delta$ is
nondegenerate and $j_A\times j_G{\vert}$ is faithful, we say \emph{Mansfield
imprimitivity works} for $N$ and $\delta$.

In Section~4 we show that Mansfield imprimitivity passes to twisted
crossed products, extending the above-mentioned result of Phillips and
Raeburn.  More explicitly, if $(A,G,G/K,\delta,\tau)$ is a twisted
coaction and $N$ is a closed normal subgroup of $G$ contained in $K$
such that Mansfield imprimitivity works for $N$ and $\delta$, then a
quotient of $Y_{G/N}^G$ is an $A\times_{G/K}G\times_rN$ --
$A\times_{G/K}G/N$ imprimitivity bimodule (\thmref{MPR-imp-thm}).

We were led to these ideas by an investigation \cite{KQR-DI}\ of the
duality of induction and restriction for coactions, and hence by a need
for a general, workable imprimitivity framework for Mansfield
induction.  In Section~5 we prove several results concerning the
compatibility of Mansfield imprimitivity with coaction constructions
such as Morita equivalence, inflation, and stabilization.  These
results will be needed in \cite{KQR-DI}; for the present, they
serve to illustrate
the robustness of Mansfield imprimitivity.

This research was carried out while the second author was visiting the
University of Newcastle in 1994 and 1995, and while the first author
was visiting Arizona State University in June, 1995. The authors thank
their respective hosts for their hospitality; the second author
particularly acknowledges Iain Raeburn. The authors further thank
Professor Raeburn for many helpful conversations.


\section{Preliminaries}

Throughout, $G$ will be a locally compact group with modular function
$\Delta_G$ and left Haar measure $ds$. The group \cstar algebra of $G$
is denoted $C^*(G)$; a subscript $r$, as in $C^*_r(G)$ or $B\times_r G$,
always indicates a reduced object.


\subsection*{Coactions}

Our coactions use the conventions of \cite{QuiFR}, \cite{QR-IC}, and
\cite{RaeCR}, although the latter uses maximal tensor products.
A (full) \emph{coaction} of $G$ on $A$ is an injective,
nondegenerate homomorphism $\delta$ from $A$ to $M(A\otimes C^*(G))$
(where we use the minimal tensor product for \cstar algebras throughout)
such that:
\begin{enumerate}
\item $\delta(a)(1\otimes z),(1\otimes z)\delta(a)\in A\otimes C^*(G)$
for all $a\in A,z\in C^*(G)$;
\item
$(\delta\otimes\iota)\circ\delta=(\iota\otimes\delta_G)\circ\delta$,
\end{enumerate}
where $\delta_G\colon C^*(G)\to M(C^*(G)\otimes C^*(G))$ is the
integrated form of the representation $s\mapsto s\otimes s$ of $G$.
The coaction $\delta$ is \emph{nondegenerate} if
$\overline{\delta_{A(G)}(A)} = A$, where
$\delta_f(a) = (\iota\otimes f)(\delta(a))$ for $f$ in the Fourier
algebra $A(G)$ and $a$ in $A$.  Equivalently, $\delta$ is nondegenerate
if $\overline{\delta(A)(1\otimes C^*(G))} = A\otimes C^*(G)$
\cite[Theorem~5]{KatTD}.

Suppose $K$ is a closed normal subgroup of $G$, and $(A,K,\epsilon)$ is
a coaction.  Let $i_K\colon C^*(K)\to M(C^*(G))$ be the canonical
nondegenerate homomorphism.  Then by \cite[Example~2.4]{PR-CO},
$$\inf\epsilon = (\iota\otimes i_K)\circ\epsilon\colon A\to M(A\otimes
C^*(G))$$
is a coaction of $G$ on $A$, called the \emph{inflation} of $\epsilon$.
It turns out that inflation respects nondegeneracy of coactions:

\begin{prop}\label{inflation-nondegen-prop}
Let $K$ be a closed normal subgroup of $G$, and let $(A,K,\epsilon)$ be
a coaction.  Then $\epsilon$ is nondegenerate if and only if the
inflated coaction $(A,G,\inf\epsilon)$ is.
\end{prop}

\begin{proof}
Suppose first that $\epsilon$ is
nondegenerate.  Then
\begin{eqnarray*}
\overline{\inf\epsilon(A)(1\otimes C^*(G))}
& = & \overline{(\iota\otimes i_K)(\epsilon(A))\bigl(1\otimes
i_K(C^*(K))C^*(G)\bigr)}\\
& = & \overline{(\iota\otimes i_K)\bigl(\epsilon(A)(1\otimes
C^*(K))\bigr) (1\otimes C^*(G))}\\
& = & \overline{(A\otimes i_K(C^*(K)))(1\otimes C^*(G))}\\
& = & A\otimes C^*(G),
\end{eqnarray*}
so $\inf\epsilon$ is nondegenerate.

Conversely, suppose $\inf\epsilon$ is nondegenerate.  A simple
calculation shows that $(\inf\epsilon)_f = \epsilon_{f{\vert}_K}$ for $f\in
A(G)$.  Now $\{ f{\vert}_K \mid f\in A(G) \}\subset A(K)$, since $A(G)$ is
the closure in $B(G)$ of $B(G)\cap C_c(G)$.  Hence,
$$A = \overline{(\inf\epsilon)_{A(G)}(A)} \subset
\overline{\epsilon_{A(K)}(A)},$$
so $\epsilon$ is nondegenerate.
\end{proof}

A \emph{covariant representation } of
a coaction $(A,G,\delta)$ is a pair
$(\pi,\mu)$, where $\pi$ and $\mu$ are nondegenerate representations of
$A$ and $C_0(G)$, respectively,  on Hilbert space, such that
$$(\pi\otimes\iota)\circ\delta(a)=
\ad\mu\otimes\iota(w_G)(\pi(a)\otimes 1)\righttext{for }a\in A,$$
where $w_G\in M(C_0(G)\otimes C^*(G))$ is the unitary element determined
by the canonical embedding of $G$ in $M(C^*(G))$. A \emph{crossed
product} for $(A,G,\delta)$ is a triple $(A\times_\delta G,j_A,j_G)$
such that $A\times_\delta G$ is a \cstar algebra (which we will denote by
$A\times G$ if $\delta$ is understood) and $(j_A,j_G)$ is a pair of
nondegenerate homomorphisms of $A$ and $C_0(G)$, respectively, to
$M(A\times G)$ satisfying:
\begin{enumerate}
\item for every nondegenerate representation $\rho$ of $A\times G$,
$(\rho\circ j_A,\rho\circ j_G)$ is a covariant representation of
$(A,G,\delta)$;
\item for every covariant representation $(\pi,\mu)$ of $(A,G,\delta)$
there is a representation $\pi\times\mu$ of $A\times G$ such that
$(\pi\times\mu)\circ j_A=\pi$ and $(\pi\times\mu)\circ j_G=\mu$;
\item $A\times_\delta G$ is the closed span of the
products $j_A(a)j_G(f)$ for $a\in A$, $f\in C_0(G)$.
\end{enumerate}
We write $j_A^G$ for $j_A$ and $j_G^A$ for $j_G$ when confusion due to
the presence of several coactions or groups is likely to arise.
There is an action $\hat\delta$, called the \emph{dual action}, of $G$ on
$A\times G$ such that $\hat\delta_s(j_A(a)j_G(f))=j_A(a)j_G(s\cdot f)$,
where $(s\cdot f)(t)=f(ts)$.

It turns out that, for any nondegenerate representation $\pi$ of $A$,
$((\pi\otimes\lambda)\circ\delta,1\otimes M)$, where $\lambda$ is the
left regular representation of $G$ and $M$ is the multiplication
representation of $C_0(G)$ on $L^2(G)$, is a covariant representation of
$(A,G,\delta)$, which moreover can be taken to be $(j_A,j_G)$ whenever
$\ker\pi\subset\ker j_A$. We call $\delta$ \emph{normal} if $j_A$ is
faithful. If $(\pi,\mu)$ is a covariant representation, then
$\ad\mu\otimes\iota(w_G)\circ(\cdot\otimes 1)$ is a normal coaction on
$\pi(A)$. In particular, $\ad j_G\otimes\iota(w_G)\circ(\cdot\otimes 1)$
is a normal coaction on $j_A(A)$ with (essentially) the same covariant
representations and crossed product as $\delta$, and we call this
coaction on $j_A(A)$ the \emph{normalization} of $\delta$. If $\delta$
is a normal coaction, then $(\iota\otimes\lambda)\circ\delta$ is a
reduced coaction on $A$, which is nondegenerate if and only if $\delta$
is. In any event, whether $\delta$ is normal or not,
$(\iota\otimes\lambda)\circ\delta$ factors to give
a reduced coaction on $A/\ker j_A$, called the \emph{reduction} of
$\delta$. Moreover, every nondegenerate reduced coaction is the
reduction of a unique normal coaction.

For example, $(C^*(G),G,\delta_G)$ is a coaction which is normal if and
only if $G$ is amenable. $(\K(L^2(G)),\lambda,M)$ is a crossed product,
so $\lambda_s\mapsto\lambda_s\otimes s$ is the normalization and
$\lambda_s\mapsto\lambda_s\otimes\lambda_s$ the reduction.


\subsection*{Hilbert modules}

Our main references for Hilbert modules and Morita equivalence are
\cite{LanHC}, \cite{RieIR}, and \cite{RieUR}.
All our Hilbert modules
(except multiplier bimodules)
will be full; i.e., the closed span of the inner
products generates the $C^*$-algebra.

\begin{defn}
A \emph{right-Hilbert $A$ -- $B$
bimodule} (a term coined by Bui in \cite{BuiFC})
is a (right) Hilbert $B$-module $X$ together with a
nondegenerate action of $A$ by adjointable $B$-module maps (i.e., there
is a homomorphism of $A$ into $\L_B(X)$ such that $AX=X$).
\end{defn}

If $X$ is also a left Hilbert $A$-module
(in the obvious sense)
such that $\lip{A}{x}{y} z = x
\rip{B}{y}{z}$ for $x,y,z\in X$, then of course $X$ is an $A$ -- $B$
imprimitivity bimodule.  We write ${}_AX_B$ when we want to emphasize
(or merely indicate) the coefficient algebras.  We denote the reverse
bimodule by $\tilde{X}$, with elements $\tilde{x}$.

For this work (and in \cite{KQR-DI}), we feel that right-Hilbert
bimodules are the right objects to use.
For example, in applications of the Rieffel induction
process it is really right-Hilbert bimodules that are required.
Routine calculations suffice to adapt most results about
Hilbert modules or imprimitivity bimodules to the setting of
right-Hilbert bimodules, and we will use such adapted results
with only a reference to the original results.

Following \cite{ER-MI},
a \emph{multiplier} $m=(m_A,m_B)$ of an $A$ -- $B$
imprimitivity bimodule $X$
consists of an $A$-linear map
$m_A\colon A\to X$ and a $B$-linear map
$m_B\colon B\to X$
such that
$ m_A(a)\cdot b = a\cdot m_B(b)$
for $a\in A$ and $b\in B$.
The \emph{multiplier bimodule} $M(X)$ consists of all multipliers of
$X$; it is naturally a Hilbert $M(A)$ -- $M(B)$ bimodule, but is in
general \emph{not} full (cf. \cite[\S1]{ER-MI}).  Following
\cite[Definition~1.8]{ER-MI}\ and \cite[Definition~A1(b)]{NgCC}, an
\emph{imprimitivity bimodule homomorphism}
$\phi=(\phi_A,\phi_X,\phi_B)\colon
{}_AX_B\to M({}_CY_D)$
consists of homomorphisms $\phi_A\colon A\to M(C)$ and
$\phi_B\colon B\to M(D)$ and a linear map $\phi_X\colon X\to M(Y)$
satisfying
\begin{itemize}
\openup1\jot
\item[(i)] $\phi_A(\lip{A}{x}{y}) = \lip{M(C)}{\phi_X(x)}{\phi_X(y)}$;
\item[(ii)] $\phi_B(\rip{B}{x}{y}) = \rip{M(D)}{\phi_X(x)}{\phi_X(y)}$;
\item[(iii)] $\phi_X(a\cdot x\cdot b) = \phi_A(a)\cdot \phi_X(x)\cdot
\phi_B(b),$
\end{itemize}
for all $a\in A$, $x,y\in X$, and $b\in B$.

Following
\cite[Definition~3.3]{NgCC}
(see also \cite{BS-CH}, \cite{BuiTC}, \cite{ER-ST}),
a \emph{coaction}\ $\delta$ of $G$ on an
imprimitivity bimodule ${}_AX_B$ is an imprimitivity bimodule
homomorphism
$$\delta =
(\delta_A,\delta_X,\delta_B)\colon {}_AX_B\to M\left({}_{A\otimes
C^*(G)}(X\otimes C^*(G))_{B\otimes C^*(G)}\right)$$
such that $(A,G,\delta_A)$ and $(B,G,\delta_B)$ are $C^*$-coactions, and
satisfying
$$(\delta_X\otimes\iota)\circ\delta_X =
(\iota\otimes\delta_G)\circ\delta_X$$
and
$$\overline{\delta_X(X)(B\otimes C^*(G))}
= X\otimes C^*(G).$$
When such a $\delta$ exists we say $(A,G,\delta_A)$ and
$(B,G,\delta_B)$ are \emph{Morita equivalent}, and we call $(X,\delta_X)$
a \emph{Morita equivalence} of $\delta_A$ and $\delta_B$.
Note that 
we automatically have
$\delta_X(x)\cdot(1_B\otimes z)$ and $(1_A\otimes z)\delta_X(x)
\in X\otimes C^*(G)$ for $x\in X$, $z\in C^*(G)$.

It turns out that Morita equivalence preserves nondegeneracy of
\cstar coactions:

\newcommand{\twoby}[4]{ \left( \begin{array}{cc}
        {#1}    &       {#2}    \\
        {#3}    &       {#4}    \end{array} \right) }

\newcommand{\smtwoby}[4]{ \bigl( \begin{smallmatrix}
{#1}&{#2}\\ {#3}&{#4}
\end{smallmatrix} \bigr) }

\begin{prop}\label{me-nondegen-prop}
Let $(A,G,\delta_A)$ and $(B,G,\delta_B)$
be Morita equivalent coactions.  Then $\delta_A$ is a 
nondegenerate coaction if
and only if $\delta_B$ is.
\end{prop}

\begin{proof}
Suppose $(X,\delta_X)$ is a Morita equivalence for
the \cstar coactions $(A,G,\delta_A)$ and
$(B,G,\delta_B)$; so $\delta = (\delta_A, \delta_X, \delta_B)$ is a
coaction of $G$ on ${}_AX_B$.

Assume $\delta_B$ is nondegenerate.  
Then $\overline{\delta_B(B)(1\otimes C^*(G))} = B\otimes C^*(G)$, so we
have
\begin{align*}
\overline{\delta_X(X)(1\otimes C^*(G))}
& = \overline{\delta_X(X\cdot B)(1\otimes C^*(G))}\\
& = \overline{\delta_X(X)\delta_B(B)(1\otimes C^*(G))}\\
& = \overline{\delta_X(X)(B\otimes C^*(G))}.
\end{align*}
Thus, 
\begin{equation}\label{star-eqn}
\overline{\delta_X(X)(1\otimes C^*(G))} = X\otimes C^*(G).
\end{equation}

We will need to slice $M(X\otimes C^*(G))$ into $M(X)$.  For this we
have found it most convenient to use the linking algebra
$$L = \twoby{A}{X}{\tilde X}{B}$$
of \cite{BGR-SI}, which by
\cite[Appendix]{ER-MI}\ carries a coaction defined by
$$\delta_L = \twoby{\delta_A}{\delta_X}{\delta_{\tilde X}}{\delta_B}.$$
Let $p=\smtwoby{1}{0}{0}{0}$ and $q=\smtwoby{0}{0}{0}{1}$.  We regard $A$,
$B$, and $X$ as sitting inside $L$, so in particular $X = pLq$.
\cite[Proposition~A.1]{ER-MI}\ shows that
$$M(L) = \twoby{M(A)}{M(X)}{M(\tilde X)}{M(B)},$$
and we have $M(X) = p M(L) q$.
Moreover,
$$L\otimes C^*(G) \cong \twoby{A\otimes C^*(G)}{X\otimes C^*(G)}{\tilde
X\otimes C^*(G)}{B\otimes C^*(G)},$$
and we blur the distinction between the two sides of this isomorphism.
Hence, we have
$$X\otimes C^*(G) = (p\otimes 1)(L\otimes C^*(G))(q\otimes 1)$$
and
$$M(X\otimes C^*(G)) = (p\otimes 1)(M(L\otimes C^*(G)))(q\otimes 1).$$

Now let $f\in A(G)$.  We can certainly use $\iota\otimes f$ to slice
$M(L\otimes C^*(G))$ into $M(L)$, and it makes sense to restrict
$\iota\otimes f$ to $M(X\otimes C^*(G))$.  We have
\begin{eqnarray*}
(\iota\otimes f)(M(X\otimes C^*(G)))
& = & (\iota\otimes f)
\bigl((p\otimes 1)M(L\otimes C^*(G))(q\otimes 1)\bigr)\\
& = & p(\iota\otimes f)(M(L\otimes C^*(G)))q\\
& \subset & p M(L) q\\
& = & M(X),
\end{eqnarray*}
so we now know how to slice $M(X\otimes C^*(G))$ into $M(X)$.

As with \cstar coactions,
we really need to slice the image of $\delta_X$ into $X$.
Without loss of
generality let $f=g\cdot d$, with $g\in A(G)$ and $d\in C^*(G)$.  Then
\begin{eqnarray*}
(\iota\otimes f)(\delta_X(X))
& = & (\iota\otimes g)((1\otimes d)\delta_X(X))\\
& \subset & (\iota\otimes g)(X\otimes C^*(G))\\
& = & (\iota\otimes g)
\bigl((p\otimes 1)(L\otimes C^*(G))(q\otimes 1)\bigr)\\
& = & p(\iota\otimes g)(L\otimes C^*(G))q\\
& \subset & pLq\\
& = & X.
\end{eqnarray*}
Hence, $X$ is an $A(G)$-submodule of $L$, as desired.

For $x\in (\iota\otimes A_c(G))(\delta_X(X))$, $f\in A_c(G)$, and $d\in
C^*(G)$, \cite[Equation~(1.2)]{QuiFR}\ tells us that
$$x\otimes u_G(f)d
= \int \delta_X\bigl((\iota\otimes s\cdot f)(x)\bigr)(1\otimes
sd)\, ds,$$
where $u_G\colon L^1(G)\to C^*(G)$ is the canonical embedding, $(s\cdot
f)(t) = f(ts)$, and the integral is norm-convergent.  Thus, the proof
of \cite[Corollary~1.5]{QuiFR}\ shows that \eqref{star-eqn}\ implies
\begin{equation}\label{two-star-eqn}
\overline{(\iota\otimes A(G))(\delta_X(X))} = X.
\end{equation}

Now, \eqref{two-star-eqn}\ is a symmetric condition which is actually
equivalent to the asymmetric \eqref{star-eqn}, so we can expect to
derive a mirror-image of \eqref{star-eqn}\ from it.
Indeed, a calculation similar to the proof of
\cite[Equation~(1.2)]{QuiFR}\ shows that
$$x\otimes du_G(f)
= \int (1\otimes ds)\delta_X\bigl((\iota\otimes(f\cdot
s)')(x)\bigr)\, ds,$$
where $(f\cdot s)'(t) = f(st)\Delta(t)$.  Thus, the proof of
\cite[Corollary~1.4]{QuiFR}\ shows that \eqref{two-star-eqn}\ implies
\begin{equation}\label{three-star-eqn}
\overline{(1\otimes C^*(G))\delta_X(X)} = X\otimes C^*(G).
\end{equation}
(Note that the proof of \cite[Corollary~1.4]{QuiFR}\ did not need to
use \cite[Equation~(1.1)]{QuiFR}.)
Using \eqref{three-star-eqn}\ in the $M(A\otimes C^*(G))$-valued inner
product on $X\otimes C^*(G)$ gives
$$\overline{(1\otimes C^*(G))\delta_A(A)(1\otimes C^*(G))} = A\otimes
C^*(G),$$
which implies $\delta_A$ is nondegenerate by an argument similar to the
proof of \cite[Corollary~1.5]{QuiFR}.

We have shown that nondegneracy of $\delta_B$ implies nondegeneracy of
$\delta_A$; by symmetry, this completes the proof.
\end{proof}


\section{Mansfield Imprimitivity}

Let $(A,G,\delta)$ be a coaction. Further let $N$ be a closed normal
subgroup of $G$, $q_N\colon C^*(G)\to C^*(G/N)$ the canonical quotient
homomorphism,
and $p_N\colon C_0(G/N)\to M(C_0(G))$ the nondegenerate embedding
$p_N(f)(s) = f(sN)$.  (When confusion is unlikely, we suppress $p_N$
and identify $C_0(G/N)$ with its image in $C_b(G)$.)
Then $\delta{\vert}:=(\iota\otimes q_N)\circ\delta$ is a coaction
of $G/N$ on $A$, called the \emph{restriction} of $\delta$ to $G/N$.
Only injectivity is nonobvious, but for this note that
$(\delta\otimes\iota)\circ\delta{\vert}
=(\iota\otimes\delta_G{\vert})\circ\delta$,
and $\delta_G{\vert}$ is injective since it has left inverse
$\iota\otimes\pi_0$, where $\pi_0\colon G/N\to\{1\}$ is the trivial
character. If $\mu$ is a homomorphism of $C_0(G)$, we
write  $\mu{\vert}$
for $\mu\circ p_N$ when confusion seems unlikely, and think of $\mu{\vert}$
as
the restriction of $\mu$ to $C_0(G/N)$.
An easy computation, using the
identity $(\iota\otimes q_N)(w_G)=
(p_N\otimes\iota)(w_{G/N})$, shows that if $(\pi,\mu)$ is a
covariant representation of $(A,G,\delta)$, then $(\pi,\mu{\vert})$ is a
covariant representation of $(A,G/N,\delta{\vert})$. In particular, since
$C_0(G/N)$ sits nondegenerately inside $M(C_0(G))$,
$$j_A\times j_G{\vert}\colon A\times G/N \to M(A\times G)$$
is a nondegenerate homomorphism.
However, $j_A\times j_G{\vert}$ will be unfaithful in general. For example,
if $\delta$ is any nonnormal coaction, then $A\times G/G=A$ and
$j_A\times j_G{\vert}=j_A$ is unfaithful.

\begin{lem}
$j_A\times j_G{\vert}\colon A\times G/N \to M(A\times G)$ is faithful
if and only if
$\ker j^G_A\subset\ker j^{G/N}_A$,
and in this case we actually have
$\ker j^G_A=\ker j^{G/N}_A$.
\end{lem}

\begin{proof}
Since $(j_A^G,j_G{\vert})$ is a covariant representation of
$(A,G/N,\delta{\vert})$, we always have
$\ker j_A^{G/N}\subset \ker j_A^G$, giving the second
statement.

The forward implication of the first statement is clear:
if $j_A\times j_G{\vert}$ is faithful, certainly
$\ker j_A^G \subset \ker j_A^{G/N}$.
For the reverse implication, assume $\ker
j^G_A\subset\ker j^{G/N}_A$. Since the dual action $\hat\delta$
leaves the image of $j_A\times j_G{\vert}$ invariant, and since the
restriction $\hat\delta{\vert}_N$ is trivial on $\im j_A\times j_G{\vert}$, we
get
an action of $G/N$ on $\im j_A\times j_G{\vert}$ such that $j_A\times
j_G{\vert}$ is
$G/N$-equivariant. Since $\ker j^G_A\subset\ker j^{G/N}_A$,
\cite[Proposition 3.9]{QuiFR} shows
$j_A\times j_G{\vert}$ is faithful.
\end{proof}

There are certainly lots of situations where $j_A\times j_G{\vert}$ is
faithful:

\begin{lem}\label{normal-or-amenable-lem}
$j_A\times j_G{\vert}\colon A\times G/N\to M(A\times G)$
is faithful if either $\delta$ is normal,
in which case $\delta{\vert}$ is normal as well, or $N$ is amenable.
\end{lem}

\begin{proof}
If $\delta$ is normal, then $\ker j^G_A=\{0\}$, so by the preceding
lemma $\ker j^{G/N}_A=\{0\}$ as well, hence $j_A\times j_G{\vert}$ is
faithful and $\delta{\vert}$ is normal.

On the other hand, if $N$ is amenable, then
$\ker\lambda_G\subset\ker\lambda_{G/N}\circ q_N$, so
\begin{align*}
\begin{split}
\ker j^G_A
&=\ker(\iota\otimes\lambda_G)\circ\delta
\subset\ker(\iota\otimes\lambda_{G/N})\circ(\iota\otimes q_N)\circ\delta\\
&=\ker(\iota\otimes\lambda_{G/N})\circ\delta{\vert}
=\ker j^{G/N}_A,
\end{split}
\end{align*}
so again using the preceding lemma, $j_A\times j_G{\vert}$ is faithful.
\end{proof}

We now adapt Mansfield's imprimitivity to our setting. Let
$(A,G,\delta)$ be a coaction and $N$ a closed normal subgroup of $G$. Since
$(j_A,j_G)$ can be taken to be
$((\iota\otimes\lambda)\circ\delta,1\otimes M)$, where $A$ is faithfully
and nondegenerately represented on a Hilbert space $\H$, all of
Mansfield's computations with $A\times G$ on $\H\otimes L^2(G)$ can be
carried out abstractly with $(j_A,j_G)$. In Mansfield's setting (with
$\delta$ reduced and $N$ amenable),
$((\iota\otimes\lambda)\circ\delta)\times(1\otimes M{\vert})$ is a faithful
representation of $A\times G/N$ \cite[Proposition 7]{ManJF}. In our
setting, this representation of $A\times G/N$ need not be faithful.
However, all of Mansfield's computations are really done with the image
of $A\times G/N$ in $B(\H\otimes L^2(G))$. In our abstract setting, we use
$\im j_A\times j_G{\vert}$ (where $j_A$ means $j^G_A$ by default), and
again all of Manfield's computations carry over.

Mansfield develops the bulk of his imprimitivity machinery in \S3 of
\cite{ManJF};
note that he does not require $N$ to be amenable for this
section. However, he implicitly requires $\delta$ to be nondegenerate
(and this becomes explicit in his proof of \cite[Theorem 12]{ManJF}),
so we will require this also from now on.

Let $A_c(G)=A(G)\cap C_c(G)$, and for a compact subset $E$ of $G$ let
$C_E(G)=\{f\in C(G)\mid\supp f\subset E\}$. For fixed $u\in A_c(G)$ and
compact $E\subset G$, say an element of $M(A\times G)$ is $(u,E,N)$ if
it is in the closed span of the products $j_A(\delta_u(a))j_G(\phi(f))$
for $a\in A$, $f\in C_E(G)$, where $\delta_u=(\iota\otimes
u)\circ\delta\colon A\to A$ and $\phi\colon C_c(G)\to C_c(G/N)$ is the
surjection
$$\phi(f)(sN)=\int_N f(sn)\,dn.$$
For $N=\{e\}$ just say the element is $(u,E)$. Manfield's computations
show the set $\D_N$ of elements of $M(A\times G)$ which are $(u,E,N)$
for some $u\in A_c(G)$ and compact $E\subset G$ is a dense \star
subalgebra of $\im j_A\times j_G{\vert}$; write $\D$ for
$\D_{\{e\}}$. We emphasize that when $A\times G$ is represented on
$\H\otimes L^2(G)$ via $(\iota\otimes\lambda)\circ\delta\times(1\otimes
M)$ we get exactly Mansfield's $\D$ and $\D_N$. Mansfield's computations
also show there is a linear map $\Psi\colon\D\to\D_N$ such that
$$\Psi(j_A(a)j_G(f))=j_A(a)j_G(\phi(f))
\righttext{for }a\in\delta_{A_c(G)}(A),f\in C_c(G).$$
Moreover, for $x,y\in \D$ the maps $n\mapsto\hat\delta_n(x)y$ and
$n\mapsto x\hat\delta_n(y)$ are in $C_c(N,\D)$, and
\begin{align*}
\Psi(x)y&=\int_N \hat\delta_n(x)y\,dn\\
x\Psi(y)&=\int_N x\hat\delta_n(y)\,dn.
\end{align*}
$\D$ becomes a (full) pre-Hilbert $\D_N$-module under right
multiplication (which makes sense since $\D_N\subset M(A\times G)$) and
inner product
$$\langle x,y\rangle_{\D_N}=\Psi(x^*y).$$
$\D$ becomes a left $C_c(N,\D)$-module under the integrated form of the
left $\D$ multiplication and the $N$-action
$$n\cdot x=\Delta(n)^{\frac12}\hat\delta_n(x)
\righttext{for }n\in N,x\in\D.$$

Let $Y^G_{G/N}$ denote the completion of the pre-Hilbert $\D_N$-module
$\D$. Then $Y^G_{G/N}$ is a right-Hilbert $A\times G\times N$ -- $\im
j_A\times j_G{\vert}$ bimodule. The computations of
\cite[proof of Proposition 26]{ManJF}
show that the homomorphism of $A\times G\times
N$ to $\L_{\im j_A\times j_G{\vert}}(Y^G_{G/N})$ has image $\K(Y^G_{G/N})$ and
the same kernel as the
regular representation $A\times G\times N\to A\times G\times_r
N$. He accomplishes the latter by showing that there is a faithful
representation of $\K(Y^G_{G/N})$, induced from the identity
representation of $\im j_A\times j_G{\vert}$ (when $(j_A,j_G)$ is taken to be
$((\iota\otimes\lambda)\circ\delta,1\times M)$), such that the
composition with $A\times G\times N\to\K(Y^G_{G/N})$ is equivalent to
the regular representation. Therefore, in our setting Mansfield's
imprimitivity theorem becomes:

\begin{thm}
\cite[Theorem 27]{ManJF} If $(A,G,\delta)$ is a nondegenerate
coaction and $N$ is a closed normal subgroup of $G$, then $Y^G_{G/N}$ is an
$A\times G\times_r N$--$\im j_A\times j_G{\vert}$ imprimitivity bimodule.
\end{thm}

We would really like a Morita equivalence between $A\times G\times_r N$
and $A\times G/N$ itself, so clearly we need exactly the condition that
$j_A\times j_G{\vert}\colon A\times G/N\to M(A\times G)$ is faithful.

\begin{cor}\label{m-imprim-cor}
Let $(A,G,\delta)$ be a nondegenerate coaction and $N$ a closed normal
subgroup of $G$ such that
$$j_A\times j_G{\vert}\colon A\times G/N\to M(A\times G)$$
is faithful.  Then $Y^G_{G/N}$ is an $A\times G\times_r
N$--$A\times G/N$ imprimitivity bimodule.
\end{cor}

In view of \corref{m-imprim-cor}, we make the following
definition:

\begin{defn}\label{M-works-def}
If $(A,G,\delta)$ is a nondegenerate
coaction and $N$ is a closed normal
subgroup of $G$, we say
\emph{Mansfield imprimitivity works} for $N$ and $\delta$ if
$j_A\times j_G{\vert}$ is faithful.
\end{defn}

When Mansfield imprimitivity works we let
$\langle\cdot,\cdot\rangle_{A\times G/N}$ denote the extension to
$Y^G_{G/N}$ of the inner product $\langle\cdot,\cdot\rangle_{\D_N}$ on
$\D$. Mansfield's computations show that the left inner product
${}_{A\times G\times_r N}\langle x,y\rangle$ for $x,y\in\D$ can be
identified with the element
$${}_{A\times G\times_r N}\langle x,y\rangle(n)
= x\hat\delta_n(y^*)\, \Delta(n)^{-\frac12}$$
of $C_c(N,\D)$.

The following lemma shows the strong connection between Mansfield
imprimitivity working and the existence of a twist (see Section~4):

\begin{lem}\label{k-works-d-normal-lem}
If $(A,G,G/K,\delta,\tau)$ is a nondegenerate
twisted coaction, then Mansfield
imprimitivity works for $K$ and $\delta$
if and only if $\delta$ is normal.
\end{lem}

\begin{proof}
Since $(A,G/K,\delta{\vert})$ is unitary, it is normal, so by definition
$$j_A^{G/K}\colon A\to M(A\times G/K)$$
is faithful.  If Mansfield imprimitivity works for $K$, then
$$j_A^G\times j_G{\vert} \colon A\times G/K\to M(A\times G)$$
is faithful.  Hence, the composition
$$j_A^G = (j_A^G\times j_G{\vert})\circ j_A^{G/K}$$
is faithful as well, so $\delta$ is normal.

The converse is part of \lemref{normal-or-amenable-lem}.
\end{proof}

To round out this discussion of Mansfield imprimitivity, we mention that
we have been unable to find an example where $\delta$ is nonnormal and
$N$ is nonamenable, but Mansfield imprimitivity still works.


\section{Imprimitivity for Twisted Coactions}

We now show how Mansfield imprimitivity passes to twisted crossed
products,
extending \cite[Theorem 4.1]{PR-CO}.
If $K$ is a closed normal subgroup of $G$, a coaction $(A,G,\delta)$
is \emph{twisted} over $G/K$ if there is a nondegenerate homomorphism
$\tau\colon C_0(G/K)\to M(A)$, called the \emph{twist}, such that:
\begin{enumerate}
\item $\delta{\vert}=\ad \tau\otimes\iota(w_{G/K})\circ(\cdot\otimes 1)$;
\item $\delta\circ\tau=\tau\otimes 1$;
\end{enumerate}
alternatively we call $(A,G,G/K,\delta,\tau)$ a \emph{twisted
coaction}. A twisted coaction $(A,G,G/K,\delta,\tau)$ is
\emph{nondegenerate}
or \emph{normal} if the untwisted coaction
$(A,G,\delta)$ is. A covariant
representation $(\pi,\mu)$ of $(A,G,\delta)$ \emph{preserves the twist}
if $\pi\circ\tau=\mu{\vert}$. The \emph{twisting ideal} of $A\times G$ is
\begin{align*}
\begin{split}
I_\tau=\bigcap&\{\ker \pi\times\mu\mid
\text{$(\pi,\mu)$ is a covariant representation}\\
&\qquad\text{of $(A,G,\delta)$ preserving the twist}\}.
\end{split}
\end{align*}
The \emph{twisted crossed product} is
$$A\times_{G/K}G=(A\times G)/I_\tau.$$
The quotient map $A\times G\to A\times_{G/K}G$ is of the form $k_A\times
k_G$ for a unique covariant representation $(k_A,k_G)$. Moreover,
$(k_A,k_G)$ preserves the twist, and for any covariant representation
$(\pi,\mu)$ preserving the twist there is a unique representation
$\pi\times_{G/K}\mu$ of $A\times_{G/K}G$ such that
$$(\pi\times_{G/K}\mu)\circ k_A=\pi
\spacetext{and}
(\pi\times_{G/K}\mu)\circ k_G=\mu.$$
The restriction $\hat\delta{\vert}_K$ of the dual action $\hat\delta$ leaves
$I_\tau$ invariant, so we get a dual action, denoted $\tilde\delta$,
of $K$ (\emph{not} $G$) on $A\times_{G/K}G$.

If $N$ is another closed normal subgroup of $G$ contained in $K$, then
$\tau$ is also a twist for the restricted coaction $(A,G/N,\delta{\vert})$,
over the quotient $G/K\cong (G/N)/(K/N)$, so there is a restricted
twisted crossed product $A\times_{G/K}G/N$. We have
$A\times_{G/K}G/K\cong A$.  When Mansfield imprimitivity works for $N$,
we want Mansfield's Hilbert module $Y^G_{G/N}$ to pass to an
$A\times_{G/K}G\times_r N$--$A\times_{G/K}G/N$ imprimitivity bimodule.
When $N=K$ is amenable (and $\delta$ is reduced), \cite{PR-CO} shows
how to do this: if $I_\tau$ and $I_\tau^K$ are the twisting ideals of
$A\times G$ and $A\times G/K$, respectively, then \cite{PR-CO} shows
$I_\tau\times K=Y^G_{G/K}\dashind I_\tau^K$, so
$Y^G_{G/K}/(Y^G_{G/K}\cdot I_\tau^K)$ gives a Morita equivalence
between $(A\times G\times K)/(I_\tau\times K)\cong A\times_{G/K}G\times
K$ and $(A\times G/K)/I_\tau^K\cong A$.  In the absence of amenability,
we run into trouble since we need to identify a quotient of the
\emph{reduced} crossed product, so we cannot use universal properties
as in \cite{PR-CO}.  We need to know the ideals of $A\times G\times_r
N$ and $A\times G/N$ match up suitably. More precisely, if $I_\tau^N$
is the twisting ideal of $A\times G/N$, we need to know
\begin{equation}\label{condition1}
(A\times G\times_r N)/\left(Y^G_{G/N}\dashind I_\tau^N\right)\cong
A\times_{G/K}G\times_r N.
\end{equation}

Let $X$ be Green's bimodule for inducing representations
from $A\times G$ to $A\times G\times N$.  $X$ is an $A\times G\times
N\times N$ -- $A\times G$ imprimitivity bimodule, which we view as a
right-Hilbert $A\times G\times_r N$ -- $A\times G$ bimodule via the
canonical map $j_{A\times G\times_r N}\colon A\times G\times_r N\to
M(A\times G\times_r N\times N) \cong M(A\times G\times N\times N)$.  By
the commutativity of the diagram
\begin{equation*}
\begin{diagram}
\node{A\times G\times N}
	\arrow{s,l}{q}
	\arrow{e,t}{j_{A\times G\times N}}
\node{M(A\times G\times N\times N)}
	\arrow{s,r}{\cong}\\
\node{A\times G\times_r N}
	\arrow{e,b}{j_{A\times G\times_r N}}
\node{M(A\times G\times_r N\times N),}
\end{diagram}
\end{equation*}
the left action of $A\times G\times_r N$ on $X$ is determined by the
left actions of $A\times G$ and $N$; we will use this when we compute
with $X$ in the proof of \thmref{special-res-ind-thm}.

Now, since
$$A\times_{G/K}G\times_r N\cong
\im \left(X\dashind (k_A^G\times k_G)\right),$$
\eqref{condition1} is equivalent to
\begin{equation}\label{condition2}
Y^G_{G/N}\dashind I_\tau^N=X\dashind I_\tau,
\end{equation}
where $I_\tau$ is the twisting ideal of $A\times G$. Since $Y^G_{G/N}$
is an $A\times G\times_r N$-$A\times G/N$ imprimitivity bimodule,
\eqref{condition2} is equivalent to
\begin{equation*}\label{condition3}
I_\tau^N=\widetilde{Y^G_{G/N}}\dashind
\left(X\dashind I_\tau\right),
\end{equation*}
which is the same as
\begin{equation}\label{condition3.1}
\ker k_A^{G/N}\times k_G = \ker \widetilde{Y_{G/N}^G}\dashind\left(
X\dashind( k_A^G\times k_G) \right).
\end{equation}
Assuming for the moment the result of \thmref{special-res-ind-thm},
this reduces to
\begin{equation}\label{condition3.2}
\ker k_A^{G/N}\times k_G = \ker k_A^G\times k_G{\vert},
\end{equation}
which becomes
\begin{equation*}\label{condition3.3}
\ker k_A^G\times_{G/K}k_G{\vert} = \{0\}
\end{equation*}
on passing to representations of the restricted twisted crossed product
$A\times_{G/K}G/N$.  Hence, in order to ensure that Mansfield
imprimitivity passes to the twisted crossed products, we need only
the fidelity of
\begin{equation}\label{condition3.4}
k_A^G\times_{G/K}k_G{\vert}\colon A\times_{G/K}G/N\to M(A\times_{G/K}G).
\end{equation}

Before we investigate this condition further, we should pause to
justify the passage from \eqref{condition3.1}\ to \eqref{condition3.2}.
This is a consequence of the following result, which is a special case
of \cite[Theorem~3.1]{KQR-DI}.
Although the proof of that particular result doesn't depend on
the results in the present paper, for the reader's peace of mind we
give here a complete proof of this special case.

\begin{thm}\label{special-res-ind-thm}
Let $(A,G,\delta)$ be a nondegenerate
coaction{\rm,} and let $N$ be  a closed normal
subgroup of $G$ such that Mansfield imprimitivity works for
$N$ {\rm(}which is automatic if $N$ is amenable{\rm)}.
Then for any representation $\pi\times\mu$ of $A\times G${\rm,} the
representation
$$\widetilde{Y_{G/N}^G}\dashind\left( X\dashind(\pi\times\mu) \right)$$
of $A\times G/N$ is unitarily equivalent to $\pi\times\mu{\vert}$.
\end{thm}

\begin{proof}
It is straightforward to check that the map $\pi\times\mu\mapsto
\pi\times\mu{\vert}$ of $\Rep A\times G$ into $\Rep A\times G/N$ is
implemented by viewing $A\times G$ as a right-Hilbert $A\times G/N$ --
$A\times G$ bimodule, using the map $j_A\times j_G{\vert}\colon A\times
G/N\to M(A\times G)$.  Hence, in order to establish the theorem, it
suffices to show that
$$\widetilde{Y_{G/N}^G}\otimes_{A\times G\times_r N}X \cong A\times G$$
as a right-Hilbert $A\times G/N$ -- $A\times G$ bimodule.  We will
actually prove the assertion that
$$Y_{G/N}^G\otimes_{A\times G/N} A\times G \cong X$$
as a right-Hilbert $A\times G\times_r N$ -- $A\times G$ bimodule, which
is equivalent because $Y_{G/N}^G$ is an imprimitivity bimodule.

Both bimodules in the tensor product are completions of Mansfield's
dense subalgebra $\D$ of $A\times G$ for the appropriate inner
products.  Our isomorphism will be the extension to
$Y_{G/N}^G \otimes_{A\times G/N} A\times G$ of the map
$\Phi\colon \D\odot \D \to C_c(N,\D)$ defined by
$$ \Phi(x\otimes{y})(h) = x \hat{\delta}_h(y).$$
Note that except for the modular function,
$\Phi(x\otimes{y})$ is just 
Mansfield's left $C_c(N,\D)$-valued inner product
${}_{A\times G\times_r N}\langle x,y^*\rangle$, and hence does indeed
give an element of $C_c(N,\D)$.
In fact, if we define $f'(n) = f(n)\Delta_N(n)^{-\frac12}$, 
then the map $f\mapsto f'$ is a homeomorphism of
$C_c(N,A\times G)$ (with the inductive limit topology) onto itself, which
takes $\Phi(\D\odot \D)$ to 
${}_{A\times G\times_r N}\langle \D,\D \rangle.$
This latter set is dense in $C_c(N,A\times G)$ for the inductive limit
topology (\cite[Lemma~25]{ManJF}); it follows that the range of $\Phi$
is also inductive limit dense in $C_c(N,A\times G)$, and therefore in $X$. 

It only remains to show that $\Phi$ 
preserves the Hilbert module structure.
For the left action of $A\times G\times_r N$,
fix $d\in \D\subset A\times G$ and $h,t\in N$.  Then:
\begin{eqnarray*}
d\cdot\Phi(x\otimes{y})(h)
 & = & d\, \Phi(x\otimes{y})(h) \\
 & = & d x \deltahat_h(y) \\
 & = & \Phi(dx\otimes{y}) \\
 & = & \Phi(d\cdot x\otimes{y})
\end{eqnarray*}
and
\begin{eqnarray*}
t\cdot\Phi(x\otimes{y})(h)
 & = & \deltahat_t\left(
\Phi(x\otimes{y})(t^{-1}h)\right) \Delta(t)^{\half} \\
 & = & \deltahat_t\left( x\deltahat_{t^{-1}h}(y)\right)
\Delta(t)^{\half} \\
 & = & \deltahat_t(x) \deltahat_h(y) \Delta(t)^{\half} \\
 & = & (t\cdot x) \deltahat_h(y) \\
 & = & \Phi( t\cdot x\otimes{y})(h).
\end{eqnarray*}
For the right action of $A\times G$,
fix $d\in \D\subset A\times G$ and $h\in N$.  Then:
\begin{eqnarray*}
\Phi(x\otimes{y})\cdot d(h)
 & = & \Phi(x\otimes{y})(h)\deltahat_h(d) \\
 & = & x \deltahat_h(yd) \\
 & = & \Phi(x\otimes(yd)(h) \\
 & = & \Phi(x\otimes{y}\cdot d)(h).
\end{eqnarray*}
For the right $A\times G$-valued inner product,
compute:
\begin{eqnarray*}
\allowdisplaybreaks
\lefteqn{\rip{A\times G}{\Phi(x\otimes{y})}{\Phi(z\otimes{w})} }\\
 & = & \int_N \hat{\delta}_h\left( \Phi(x\otimes{y})(h^{-1})^*
\Phi(z\otimes{w})(h^{-1})\right)\, dh \\
 & = & \int_N \deltahat_h\left( (x\deltahat_{h^{-1}}(y))^*
z\deltahat_{h^{-1}}(w) \right)\, dh \\
 & = & y^* \int_N \deltahat_h(x^*z)\, dh\, w \\
 & = & y^* \rip{A\times G/N}{x}{z} w \\
 & = & y^*\left( \rip{A\times G/N}{x}{z}\cdot w\right) \\
 & = & \rip{A\times G}{y}{\rip{A\times G/N}{x}{z}\cdot w} \\
 & = & \rip{A\times G}{x\otimes{y}}{z\otimes{w}}.
\end{eqnarray*}

This completes the proof of the theorem.
\end{proof}

Of course, the discussion preceding \thmref{special-res-ind-thm}\ was
predicated on Mansfield imprimitivity working for $N$ and $\delta$;
however, $k_A^G\times_{G/K}k_G{\vert}$ is certainly well-defined for any
twisted coaction.  It will turn out (\thmref{m-iff-mpr-thm}) that
for a twisted coaction
$(A,G,G/K,\delta,\tau)$ and a closed normal subgroup $N$ of $G$
contained in $K$,
this map is faithful exactly when $j_A^G\times j_G{\vert}$ is.
In order to see this, we will need the
following variation on \cite[Lemma~3.5]{QR-IC}.  Recall that for an
action $(B,K,\alpha)$ of a closed subgroup $K$ of a locally compact
group $G$, the \emph{induced algebra} $\Ind_K^G B$ is the $C^*$-algebra
of continuous maps $f\colon G\to B$ such that $f(sk) =
\alpha_{k^{-1}}(f(s))$ for all $s\in G$, $k\in K$, and such that the
map $sK\mapsto \|f(s)\|$ vanishes at infinity on $G/K$.

\begin{lem}\label{ind-hom-lem}
Let $(B,K,\alpha)$ and $(C,K,\beta)$ be actions of a closed
subgroup $K$ of a locally compact group $G$, and let $\phi\colon B\to
M(C)$
be a $K$-equivariant nondegenerate homomorphism.  Then $\phi$ is
faithful if and only if the induced homomorphism $\ind\phi\colon
\ind_K^G B\to M(\ind_K^G C)$ is, where $\ind\phi$ is defined by
$$\ind\phi(f)(s) = \phi(f(s)).$$
\end{lem}

\begin{proof}
If $\phi$ has nontrivial kernel, then $\ker\phi$ is a nonzero
$K$-invariant ideal of $B$, so $\ind(\ker\phi)$ is a nonzero subset of
$\ker(\ind\phi)$: for $f\in\ind(\ker\phi)$ and $s\in G$,
$$\ind\phi(f)(s) = \phi(f(s)) = 0.$$
Conversely, if $f$ is a nonzero element of $\ker(\ind\phi)$, then $f(G)$
is a nonzero subset of $\ker\phi$.
\end{proof}

\begin{thm}\label{m-iff-mpr-thm}
Let $(A,G,G/K,\delta,\tau)$ be a twisted coaction, and
let $N$ be a closed normal subgroup of $G$ contained in $K$.
Then
$$j_A\times j_G{\vert}\colon A\times G/N\to M(A\times G)$$
is faithful
if and only if
$$k_A\times_{G/K}k_G{\vert}\colon A\times_{G/K}G/N\to
M(A\times_{G/K}G)$$
is faithful.
\end{thm}

\begin{proof}
Note that the diagram
$$
\begin{CD}
A\times G/N @>{j_A\times j_G{\vert}}>> M(A\times G) \\
@V{k_A\times k_{G/N}}VV @VV{k_A\times k_G}V \\
A\times_{G/K}G/N @>>{k_A\times_{G/K}k_G{\vert}}> M(A\times_{G/K}G)
\end{CD}
$$
commutes.

By \cite[Theorem~4.4]{QR-IC}, the formula
$$\Phi(x)(s)=(k_A\times k_G)\circ\hat\delta_{s^{-1}}(x)
\righttext{for}x\in A\times G,s\in G$$
defines an isomorphism $\Phi\colon A\times G\to\ind_K^G
A\times_{G/K}G$,
and similarly
$$A\times G/N\cong\ind_{K/N}^{G/N}A\times_{G/K}G/N.$$
We need to do everything in terms of $K$ and
$G$ rather than their quotients by $N$. Define
$\alpha\colon G\to \aut A\times G/N$ by
$$\alpha_s=(\delta{\vert})\widehat{\ }_{sN}.$$
Then
$$(j_A\times j_G{\vert})\circ\alpha_s=
\hat\delta_s\circ(j_A\times j_G{\vert})
\righttext{for}s\in G.$$
The twisting ideal $I^N_\tau$ of $A\times G/N$ is invariant under the
restriction $\alpha{\vert}_K$, so there is a unique action $\tilde\alpha$ of
$K$ on $A\times_{G/K}G/N$ such that
$$\tilde\alpha_s\circ(k_A\times k_{G/N})=
(k_A\times k_{G/N})\circ\alpha_s
\righttext{for}s\in K.$$
Of course, $\tilde\alpha_s=({\delta{\vert}})\widetilde{\ }_{sN}$.

An easy calculation shows the formula
$$\Psi(f)(s)=f(sN)\righttext{for}s\in G$$
gives an isomorphism
$$\Psi\colon\ind_{K/N}^{G/N}A\times_{G/K}G/N\to
\ind_K^G A\times_{G/K}G/N.$$
Combining this with \cite[Theorem~4.4]{QR-IC}, we can define an
isomorphism
$$\Phi_N\colon A\times G/N\to\ind_K^G A\times_{G/K}G/N$$
by
$$\Phi_N(x)(s)=(k_A\times k_{G/N})\circ\alpha_{s^{-1}}(x).$$
The following calculation shows
$$k_A\times_{G/K}k_G{\vert}\colon A\times_{G/K}G/N\to
M(A\times_{G/K}G)$$
is $K$-equivariant: for $s\in K$
\begin{align*}
\begin{split}
\lefteqn{(k_A\times_{G/K}k_G{\vert})\circ\tilde\alpha_s\circ
(k_A\times k_{G/N})} \\
\quad&=(k_A\times_{G/K}k_G{\vert})\circ(k_A\times k_{G/N})\circ\alpha_s
\\
\quad&=(k_A\times k_G)\circ(j_A\times j_G{\vert})\circ\alpha_s \\
\quad&=(k_A\times k_G)\circ\hat\delta_s\circ(j_A\times j_G{\vert}) \\
\quad&=\tilde\delta_s\circ(k_A\times k_G)\circ(j_A\times j_G{\vert}) \\
\quad&=\tilde\delta_s\circ(k_A\times_{G/K}k_G{\vert})\circ
(k_A\times k_{G/N}).
\end{split}
\end{align*}

We next show the nondegenerate homomorphism
\begin{equation*}
\Phi\circ(j_A\times j_G{\vert})\circ\Phi^{-1}_N\colon
\ind_K^G A\times_{G/K}G/N\to M(\ind_K^G A\times_{G/K}G)
\end{equation*}
is induced from
$$k_A\times_{G/K}k_G{\vert}\colon
A\times_{G/K}G/N\to M(A\times_{G/K}G)$$
in the sense of \lemref{ind-hom-lem}, and then the lemma will yield the
present theorem. For $f\in\ind_K^G A\times_{G/K}G/N$, $s\in G$,
\begin{align*}
\begin{split}
\lefteqn{\Phi\circ(j_A\times j_G{\vert})\circ\Phi^{-1}_N(f)(s)} \\
\quad&=(k_A\times k_G)\circ\hat\delta_{s^{-1}}\circ
(j_A\times j_G{\vert})\circ\Phi^{-1}_N(f) \\
\quad&=(k_A\times k_G)\circ(j_A\times j_G{\vert})
\circ\alpha_{s^{-1}}\circ\Phi^{-1}_N(f) \\
\quad&=(k_A\times_{G/K}k_G{\vert})\circ(k_A\times k_{G/N})
\circ\alpha_{s^{-1}}\circ\Phi^{-1}_N(f) \\
\quad&=(k_A\times_{G/K}k_G{\vert})\bigl(\Phi_N(\Phi^{-1}_N(f))(s)\bigr)
\\
\quad&=(k_A\times_{G/K}k_G{\vert})(f(s)).
\end{split}
\end{align*}
\end{proof}

Now the discussion encompassing
Equations \eqref{condition1}--\eqref{condition3.4},
together with \thmref{m-iff-mpr-thm},
gives the following extension
of \cite[Theorem~4.1]{PR-TC}.
In the case where $N$ is amenable, this result is hinted at in
the discussion preceding \cite[Theorem~4.7]{ER-ST}.

\begin{thm}\label{MPR-imp-thm}
Let $(A,G,G/K,\delta,\tau)$ be a
nondegenerate
twisted coaction, and let $N$
a closed normal subgroup of $G$ contained in $K$ such that
Mansfield imprimitivity works for $N$ and $\delta$
{\rm(}which is automatic if $N$ is amenable{\rm)}.
Then the quotient
$Z^G_{G/N}=Y^G_{G/N}/(Y^G_{G/N}\cdot I_\tau^N)$ is
an $A\times_{G/K}G\times_r N$--$A\times_{G/K}G/N$ imprimitivity
bimodule.
\end{thm}


\section{Subgroups, Morita Equivalence, Inflation, Stabilization}

In this section, we will show that Mansfield imprimitivity
is compatible with many of the standard coaction constructions.  As a
warmup, we show that Mansfield imprimitivity passes to subgroups.
First we need the following
variation on \cite[Corollary~4.10]{QR-IC}.  Sadly, the
hypothesis in \cite{QR-IC}\ seems to be deficient: equivariance
must be imposed on the integrated form of the pair $(\pi,\mu)$, rather
than just $\mu$.  

\begin{prop}\label{fidelity-prop}
Let $(A,G,G/K,\delta,\tau)$ be a twisted coaction and $(\pi,\mu)$ a
covariant representation preserving the twist. Then the representation
$\pi\times_{G/K}\mu$ of $A\times_{G/K}G$ is faithful if and only if
$\ker\pi\subset\ker j_A$ and there is an action $\alpha$ of $K$ on
$\im\pi\times\mu$ such that
$$\alpha_s\circ(\pi\times_{G/K}\mu)
=(\pi\times_{G/K}\mu)\circ\tilde\delta_s
\righttext{for}s\in K.$$
\end{prop}

\begin{proof}
Replace the twisted coaction by its reduction
\cite[Corollary~3.8]{QuiFR},
\cite[Theorem~4.1]{RaeCR}; then $j_A$ is faithful, and
\cite[Corollary~4.10]{QR-IC}\ gives the proposition.
\end{proof}

\begin{thm}\label{subgroup-works-thm}
Let $(A,G,\delta)$ be a nondegenerate
coaction and $N\subset H$ closed normal subgroups of
$G$.
If Mansfield imprimitivity works for $H$, then it also works
for $N$.
\end{thm}

\begin{proof}
The diagram
$$\begin{CD}
M(C_0(G)) @<{p_N}<< M(C_0(G/N)) \\
@A{p_H}AA @AA{p_{H/N}}A \\
C_0(G/H) @<{\cong}<\gamma< C_0((G/N)/(H/N))
\end{CD}$$
commutes, where $\gamma$ is the natural isomorphism, and this gives us a
commutative diagram
$$\begin{CD}
M(A\times G) @<{j_A\times(j_G\circ p_N)}<< M(A\times G/N) \\
@A{j_A\times(j_G\circ p_H)}AA @AA{j_A\circ(j_{G/N}\circ p_{H/N})}A \\
A\times G/H @<{\cong}<{j_A\times(j_{G/H}\circ\gamma)}<
A\times(G/N)/(H/N).
\end{CD}$$
Assuming $j_A\times(j_G\circ p_H)$ is faithful, so is
$j_A\times(j_{G/N}\circ p_{H/N})$. So, \cite[Theorem~3.1]{PR-TC}\ applies,
giving a decomposition isomorphism
$$\theta\colon A\times G/N \stackrel{\cong}{\to}
A\times(G/N)/(H/N)\times_{(G/N)/(H/N)}G/N,$$
hence a homomorphism
\begin{align*}
\sigma:&=(j_A\times(j_G\circ p_N))\circ\theta^{-1} \\
&=(j_A\times(j_G\circ p_N\circ p_{H/N}))
\times_{(G/N)/(H/N)}(j_G\circ p_N)
\end{align*}
of $A\times(G/N)/(H/N)\times_{(G/N)/(H/N)}G/N$. We aim to apply the
preceding proposition to show $\sigma$ is faithful. Since the
decomposition coaction of $G/N$ on $A\times(G/N)/(H/N)$ is normal, we
must show the homomorphism $j_A\times(j_G\circ p_N\circ p_{H/N})$ of
$A\times(G/N)/(H/N)$ is faithful and there is an action $\alpha$ of
$H/N$ on $\im\sigma$ such that
$$\alpha_{hN}\circ j_G\circ p_N(f)
=j_G\circ p_N(hN\cdot f)\righttext{for}h\in H,f\in C_0(G/N).$$
The first follows from
$$j_A\times(j_G\circ p_N\circ p_{H/N})
=(j_A\times(j_G\circ p_H))\circ(j_A\times(j_{G/H}\circ\gamma))$$
and fidelity of $j_A\times(j_G\circ p_H)$. For the second, since the
restriction $\hat\delta{\vert}_N$ is trivial on
$\im\sigma=\im j_A\times(j_G\circ p_N)$,
there is a unique action $\alpha$ of $H/N$ on $\im\sigma$ such that
$$\alpha_{hN}=\hat\delta_h\righttext{for}h\in H.$$
For $f\in C_0(G/N)$,
\begin{align*}
\alpha_{hN}\circ j_G\circ p_N(f)
&=\hat\delta_h\circ j_G(p_N(f)) \\
&=j_G(h\cdot p_N(f)) \\
&=j_G\circ p_N(hN\cdot f).
\end{align*}
\end{proof}

We now show Mansfield imprimitivity is preserved by Morita equivalence
of coactions.  Our conventions are those of \cite{NgCC}.

\begin{thm}\label{me-works-thm}
Let $(A,G,\delta_A)$ and $(B,G,\delta_B)$ be Morita equivalent
coactions with one {\rm(}hence both{\rm)} nondegenerate, and let $N$ be a
closed normal subgroup of $G$. Then Mansfield
imprimitivity works for $N$ and $\delta_A$ if and only if it works for
$N$ and $\delta_B$.
\end{thm}

\begin{proof}
Let $(X,\delta_X)$ be a
Morita equivalence for $(A,G,\delta_A)$ and
$(B,G,\delta_B)$; so $\delta=(\delta_A,\delta_X,\delta_B)$ is a coaction
of $G$ on ${}_AX_B$.  Nondegeneracy of both coactions follows from
nondegeneracy of either one by \propref{me-nondegen-prop}.

As usual, let $(A\times G, j_A^G, j_G^A)$ and $(B\times G, j_B^G,
j_G^B)$ denote the crossed products for $(A,G,\delta_A)$ and
$(B,G,\delta_B)$, respectively.  By the uniqueness of the imprimitivity
bimodule crossed product (\cite[Remark~3.6(c)]{NgCC}), we can suppose
that the crossed product for $({}_AX_B, G, \delta)$ is of the form
$$({}_{A\times G}(X\times G)_{B\times G},
j_A^G, j_X^G, j_B^G, j_G^A, j_G^B)$$
for some linear map $j_X^G\colon X\to M(X\times G)$ (cf.
\cite[Definition~3.5(b)]{NgCC}).

The
conscientious reader will check that
$$(j_A^G,j_X^G,j_B^G,j_G^A{\vert},j^B_G{\vert})$$
is a covariant representation of the restricted
imprimitivity bimodule coaction $({}_AX_B,G/N,\delta{\vert})$
in the sense of \cite[Definition~3.5(a)]{NgCC}.  It then follows
that there is a unique imprimitivity
bimodule representation  
$(\phi_{A\times G/N}, \phi_{X\times G/N}, \phi_{B\times G/N})$
of
${}_{A\times G/N}(X\times G/N)_{B\times G/N}$ 
such that
\begin{align*}
\begin{split}
&(\phi_{A\times G/N}\circ j_A^{G/N}, \phi_{X\times G/N}\circ
j_X^{G/N}, \phi_{B\times G/N}\circ j_B^{G/N},\\
&\quad \phi_{A\times G/N}\circ
j_{G/N}^A, \phi_{B\times G/N}\circ j_{G/N}^B) \\
&\quad\quad =  (j_A^G, j_X^G, j_B^G, j_G^A{\vert}, j_G^B{\vert}).
\end{split}
\end{align*}
In particular, $\phi_{A\times G/N}$
is a representation of $A\times G/N$
satisfying
$$\phi_{A\times G/N}\circ j_A^{G/N} = j_A^G$$
and
$$\phi_{A\times G/N}\circ j_{G/N}^A = j_G^A{\vert},$$
so we have $\phi_{A\times G/N} = j_A^G\times j_G^A{\vert}$.  By the same
token, we have $\phi_{B\times G/N} = j_B^G\times j_G^B{\vert}$.
Now \cite[Lemma~2.7]{ER-ST}\
tells us the ideals $\ker(j_A^G\times j^A_G{\vert})$ of
$A\times G/N$ and $\ker(j_B^G\times j^B_G{\vert})$ of $B\times G/N$
correspond via $X\times G/N$; in particular, $j_A^G\times j^A_G{\vert}$ is
faithful if and only if $j_B^G\times j^B_G{\vert}$ is.  This establishes the
theorem.
\end{proof}

Next, we turn to inflation.  Recall that if $(A,G,\inf\epsilon)$ is
inflated from $(A,K,\epsilon)$, then
$\inf\epsilon$ is trivially twisted over $G/K$ by
$f\mapsto f(e)1$, and \cite[Example~2.14]{PR-TC} gives a natural
isomorphism
$$A\times_{G/K}G\cong A\times K.$$

\begin{thm}\label{inflation-works-thm}
Let $N\subset K$ be closed normal subgroups of $G$,
and let $(A,K,\epsilon)$
be a coaction such that either $\epsilon$ or $\inf\epsilon$
{\rm(}hence the other{\rm)}
is nondegenerate.
Then Mansfield imprimitivity works for $N$ and $\epsilon$ if and
only if it works for $N$ and $\inf\epsilon$.
\end{thm}

\begin{proof}
Nondegeneracy of both coactions follows from nondegeneracy of either
one by \propref{inflation-nondegen-prop}.

Consider the twisted inflated coaction
$(A,G,G/K,\inf\epsilon,1)$.
The diagram
$$\begin{CD}
A\times K/N @>{j_A\times j_K{\vert}}>> M(A\times K) \\
@V{\cong}VV @VV{\cong}V \\
A\times_{G/K}G/N @>>{k_A\times_{G/K}k_G{\vert}}>
M(A\times_{G/K}G)
\end{CD}$$
commutes, so $j_A\times j_K{\vert}$ is faithful if and only if
$k_A\times_{G/K}k_G{\vert}$ is. \thmref{m-iff-mpr-thm} tells us this latter
is equivalent to fidelity of $j_A\times j_G{\vert}\colon A\times G/N\to
M(A\times G)$, and this is enough to finish the proof.
\end{proof}

Finally, we chain the above results together to show
that Mansfield imprimitivity is
compatible with the stabilization
trick of \cite{ER-ST}.
Let $(A,G,G/K,\delta,\tau)$ be a nondegenerate
twisted coaction such that
Mansfield imprimitivity works for $\delta$ and $K$ itself. In
light of \lemref{k-works-d-normal-lem}, this is equivalent to
$\delta$ being normal.
We have a dual action $\tilde\delta$ of $K$ on the twisted
crossed product $A\times_{G/K}G$, hence a double dual coaction
$\widehat{\widetilde\delta}$ of $K$ on the full crossed product
$A\times_{G/K}G\times K$.
The normalization $(\widehat{\widetilde\delta})^{\rm n}$
of this coaction
is on the reduced crossed product
$A\times_{G/K}G\times_r K$ \cite[Proposition 3.2 (1)]{RaeCR},
\cite[Propositions 3.6 and 3.7]{QuiFR}.
By \cite[Theorem 3.1]{ER-ST}, the
twisted coaction $(A,G,G/K,\delta,\tau)$ is Morita equivalent to
the inflated twisted coaction $(A\times_{G/K}G\times_r K,G,G/K,\inf
(\widehat{\widetilde\delta})^{\rm n},1)$. Even though \cite{ER-ST} uses reduced
coactions and requires $N$ to be amenable, their arguments carry over to
our setting since we assume that Mansfield imprimitivity works. The
necessary adjustments are fairly obvious, such as replacing \cite[Lemma
3.11]{ER-ST}\ with
$$\hat\delta_n\otimes\iota(j_G\otimes\iota(w_G))
=j_G\otimes\iota(w_G)(1\otimes n)\righttext{for}n\in N.$$
We do not
need the twisting subgroup $K$ to be amenable, since we use full
coactions---see \cite[Section 7]{QR-IC}.

\begin{thm}
Let $(A,G,G/K,\delta,\tau)$ be
a nondegenerate normal twisted coaction.
Then the stabilized coaction $(A\times_{G/K}G\times_r
K,K,(\widehat{\widetilde\delta})^{\rm n})$ is also nondegenerate.
Furthermore, if $N$ is a closed normal subgroup of $G$ contained in
$K$, then Mansfield imprimitivity works for $N$ and
$(\widehat{\widetilde\delta})^{\rm n}$.
\end{thm}

\begin{proof}
By \cite[Theorem~3.1]{ER-ST},
$(A,G,\delta)$ is Morita equivalent to the inflated
coaction
$(A\times_{G/K}G\times_r K,G,\inf (\widehat{\widetilde\delta})^{\rm n})$;
thus $\inf (\widehat{\widetilde\delta})^{\rm n}$ is nondegenerate by
\propref{me-nondegen-prop}, and hence also
$(\widehat{\widetilde\delta})^{\rm n}$
by \propref{inflation-nondegen-prop}.

Now let $N$ be a closed normal subgroup of $G$ contained in $K$.
Since $\delta$ is normal,
Mansfield imprimitivity works for $N$ and $\delta$
(\lemref{normal-or-amenable-lem}), and hence it also works for
$N$ and the Morita equivalent coaction
$\inf (\widehat{\widetilde\delta})^{\rm n}$ (\thmref{me-works-thm}).
It follows that Mansfield imprimitivity works for $N$ and the deflated
coaction $(\widehat{\widetilde\delta})^{\rm n}$, by
\thmref{inflation-works-thm}.
\end{proof}


\end{document}